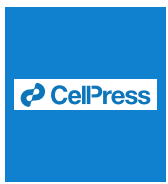
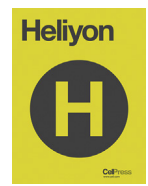

Research article

# Memory and mutualism in species sustainability: A time-fractional Lotka-Volterra model with harvesting

Mohammad M. Amirian [a,*], I.N. Towers [b], Z. Jovanoski [b], Andrew J. Irwin [a]

[a] *Department of Mathematics and Statistics, Dalhousie University, Halifax, Nova Scotia, Canada*
[b] *School of Science, UNSW Canberra, ACT, Australia*

## A R T I C L E   I N F O

## A B S T R A C T

We first present a predator-prey model for two species and then extend the model to three species where the two predator species engage in mutualistic predation. Constant effort harvesting and the impact of by-catch issue are also incorporated. Necessary sufficient conditions for the existence and stability of positive equilibrium points are examined. It is shown that harvesting is sustainable, and the memory concept of the fractional derivative damps out oscillations in the population numbers so that the system as a whole settles on an equilibrium quicker than it would with integer time derivatives. Finally, some possible physical explanations are given for the obtained results. It is shown that the stability requires the memory concept in the model.

## 1. Introduction

Ecological models are of great importance for environmental decision making because they provide stakeholders with a conceptual framework and a "laboratory" for studying the consequences of alternative policies and management scenarios [1]. One common method that scholars use to improve our understanding of environmental phenomena is the Lotka-Volterra (or predator-prey) model – an important and popular prototype model appearing in various fields of applied mathematics – due to its descriptive power, tractability and diverse applications [2, 3, 4, 5, 6, 7, 8, 9]. Therefore, many efforts have been made so far to propose more realistic models incorporating mutualism [10, 11], parasitism [12, 13], and the impact of harvesting [14, 15, 16, 17, 18, 19, 20].

Some scholars extend models formulated with fractional derivatives through applying fractional calculus (FC). In 2013, by subjecting the predators to harvesting, a modified fractional version of predator-prey model with a type II functional response was proposed [15]. Later, the fractionalised model was developed further by adding an economic interest equation to the model [21]. More recently, in addition to using the Caputo derivative to study the complex behaviour of the phenomena, scholars applied the newly formulated fractional version of the Adams-Bashforth method in their research [22, 23, 24, 25, 26, 27]. As another, for example, Kolade proposed a three-component time-fractional system representing the interaction among prey, intermediate-predators, and predators and examined the model behaviour under Caputo or the Atangana– Baleanu fractional derivatives [28].

In this study, we consider a new type of fractional model for three species (two predators and one prey) with type II mutualistic predation. An example of the type of interaction the model idealises is the mutualistic predation of spotted dolphins and yellowfin tuna upon schools of lanternfish [29]. Since the understanding of harvesting is an important issue for fishering, we assumed that all species have a market value, and also considered the impact of by-catch for dolphins when fishing for tuna. Making a comparison between the integer derivative model and fractional model under different scenarios, we examine the concept of "memory" on our model. We show that stability is more robust when the species exhibit "memory". We examine the impact of harvesting on the system with and without the memory concept as well.

## 2. Preliminaries on fractional calculus and main models

Fractional calculus is a powerful tool which has been employed in different fields of science to model complex systems with non-local be-haviour and long-term memory [30, 31, 32, 33, 34, 35]. This approach can lead to models capturing more of the phenomena under scrutiny while still keeping the model parameters to a minimum.

* Corresponding author.
*E-mail address:* M.Amirianmatlob@dal.ca (M.M. Amirian).

The adjective "fractional" in FC is a historical remnant and this calculus is a generalisation from integer order derivatives to arbitrary real-valued order, not merely rational order. Generally speaking, it can be formulated as follows.

$$\frac{d^\alpha f}{dt^\alpha} = \begin{cases} D^\alpha f(t) & \alpha > 0 \\ f(t) & \alpha = 1 \\ I^\alpha f(t) & \alpha < 0 \end{cases} \quad (1)$$

where $D^\alpha$ and $I^\alpha$ are the fractional derivative and fractional integral respectively [36]. Various definitions of fractional calculus have been proposed. All definitions coincide when the order is integer, however, this need not be the case for non-integer order. Therefore, different physical interpretations, known as the memory concept, are proposed in the fractional case [30].

We will use the Caputo definition (for other formula, see [30, 36, 37]).

$$_a^C D_t^\alpha f(t) = \int_a^t w(t-\tau) Df(\tau) d\tau \quad \text{s.t } w(t) = \frac{t^{-\alpha}}{\Gamma(1-\alpha)} \quad (2)$$

where $0 < \alpha < 1$, and $w(t)$ is a weight function whose task is the storage of the system memory over time [30].

Larger values of $\alpha$ increase the weight on the integrand $Df$ close to $t$, emphasizing the memory of nearby values of $f$. When alpha is close to zero, models formulated with the Caputo derivative ($0 < \alpha < 1$), will retain close to complete memory of the past starting at time $a$. It is expected that when the system maintains a near total memory of its past then the system resists changing over time. Based on this interpretation, therefore, we are expecting that oscillations in the population numbers of a species damp out in a system with the fractional time derivative of the order less than unity. We have examined this issue in the forthcoming sections and shown that this is the case with fractional prey-predator models.

### 2.1. Single predator model

In this section, we considered a fractional model of predator and prey. Using constant harvest quota $H_1(X) = h_1 X$ and $H_2(Y) = h_2 Y$, we assume that either both species have market value or that one species is caught as by-catch. By doing so, we tried to answer such questions as *how does the harvesting of the species affect the natural equilibrium of the ecology? How heavily can a species be harvested and still be sustainable* [20, 38]?

$$_0^C D_t^\alpha X(T) = rX\left(1 - \frac{X}{K}\right) - aXY - H_1(X),$$
$$_0^C D_t^\alpha Y(T) = \frac{aXY}{1+\sigma X} - kY - H_2(Y), \quad (3)$$

where $X$ and $Y$ are the population densities of prey and predator respectively, $T$ is time, $r$ is the prey growth rate, $K$ is the environmental carrying capacity for the prey, $a$ is the feeding rate of predators, $\sigma$ is the predator growth saturation factor and $k$ is the predator death rate. All parameters are positive reals.

After substituting the rescalings $X = Kx$, $Y = ky/a$ and $T = t/k$ in to eq. (3), we arrived at the following dimensionless system

$$\frac{d^\alpha x}{dt^\alpha} = \rho x(1-x) - xy - \varepsilon_1 x,$$
$$\frac{d^\alpha y}{dt^\alpha} = \frac{\psi xy}{1+\phi x} - y - \varepsilon_2 y, \quad (4)$$

where $\rho = r/k$, $\psi = aK/k$, $\phi = K\sigma$, $\varepsilon_1 = h_1/k$ and $\varepsilon_2 = h_2/k$.

### 2.2. Two predators model

In this section, the system (3) is extended to model interactions amongst three species: one prey and two predators. The predators are not treated as isolated hunters. Rather, we consider the predator species to be cooperative. The model of the three species with type II mutualism [10, 39] functional response for the predators is

$$_0^C D_t^\alpha X(T) = rX\left(1 - \frac{X}{K}\right) - X(aY + bZ + \xi YZ) - H_1(X),$$
$$_0^C D_t^\alpha Y(T) = \frac{XY(a+\xi Z)}{1+\sigma_1 X} - k_1 Y - H_2(Y), \quad (5)$$
$$_0^C D_t^\alpha Z(T) = \frac{XZ(b+\xi Y)}{1+\sigma_2 X} - k_2 Z - H_3(Z).$$

where $X$ is the prey population density (lanternfish), $Y$ and $Z$ are population densities of distinct predators (tuna fish and dolphin respectively), $T$ is time, $r$ is the prey growth rate, $K$ is the environmental carrying capacity for the prey, $a$, $b$, and $\xi$ are the feeding rate of predators, $\sigma_1$ and $\sigma_2$ are the predator growth saturation factor and $k_1$ and $k_2$ are the predator death rate. As before, all parameters are positive reals and the terms $H_j$ for $j = 1, 2, 3$ are the harvesting functions.

The number of parameters in system (5) can be reduced by considering the following transformations

$$X = Kx, \qquad Z = \left(\frac{a}{\xi}\right)z, \qquad Y = \frac{k_1}{a}y, \qquad T = \frac{t}{k_1}$$

thus we arrive the following dimensionless system:

$$\frac{d^\alpha x}{dt^\alpha} = \rho x(1-x) - x(y + \eta z + yz) - \varepsilon_1 x,$$
$$\frac{d^\alpha y}{dt^\alpha} = \frac{\psi xy(1+z)}{1+\phi x} - \varepsilon_2 y, \quad (6)$$
$$\frac{d^\alpha z}{dt^\alpha} = \frac{\beta xz(\eta + y)}{1+\phi_1 x} - \varepsilon_3 z.$$

where $\rho = r/k_1$, $\psi = aK/k_1$, $\beta = bK/k_1$, $\eta = ab/\xi k_1$, $\phi = \sigma_1 K$, $\phi_1 = \sigma_2 K$, $\varepsilon_1 = h_1/k_1$, $\varepsilon_2 = 1 + h_2/k_1$ and $\varepsilon_3 = (k_2 + h_3)/k_1$.

## 3. Model analysis

### 3.1. Existence and uniqueness of the non-negative solution

**Definition 3.1.** For $\alpha > 0$ and $\beta \geq 0$, the Mittag-Leffler function is defined by the following series:

$$E_{\alpha,\beta} = \sum_{k=0}^{\infty} \frac{x^k}{\Gamma(k\alpha + \beta)}$$

where $\Gamma$ is the gamma function.

**Theorem 3.1.** *Assume that* $\Omega = \{(x,y) \in \mathbb{R}_+^2 : \max\{|x|,|y|\} \leq M\}$ *and* $S = \Omega \times [t_0, T]$ *where* $T < +\infty$. *Then for any initial conditions* $(x(t_0), y(t_0)) \in \Omega$, *all the solutions* $(x(t), y(t)) \in S$ *of model* (4) *are non-negative and unique for all* $t \geq 0$.

**Proof.** (Proof by contradiction): let $k(t) = \min\{x(t), y(t)\}$, then $k(t) > 0$. Also assume

$$\exists \bar{t} > 0 \text{ s.t } k(\bar{t}) = 0, \text{ and } k(t) > 0 \quad \forall t \in [0, \bar{t}).$$

If $k(\bar{t}) = x(\bar{t})$, then after taking Laplacian transform (Table C1 & C2 of [30]) from the first equation of system (4), we have

$$x(\bar{t}) = x(0) E_\alpha \left[ (\rho(1-x) - (y+\varepsilon_1)) \right] > 0$$

This leads to a contradiction. In the similar way, when $k(\bar{t}) = y(\bar{t})$, we can obtain the contradiction. Therefore, $\forall t \geq 0$, $k(t) > 0$ and as a result, $(x(t), y(t))$ will be positive for all $t \geq 0$.

We now show that the system (4) satisfies the locally Lipschitz condition [40] needed to establish the existence and uniqueness of solutions to system (4).

Consider a mapping $\mathbb{F}(X) = (\mathbb{F}_1(X), \mathbb{F}_2(X))$ with $||.||$ norm such that

$$\mathbb{F}_1(X) = \rho x(1-x) - xy - \varepsilon_1 x,$$
$$\mathbb{F}_2(X) = \frac{\psi xy}{1+\phi x} - y - \varepsilon_2 y,$$

We show that

$$\forall X, \bar{X} \in \Omega, \ \exists L \geq 0 \ \text{s.t} \ ||\mathbb{F}(X) - \mathbb{F}(\bar{X})|| \leq L ||X - \bar{X}||$$

where $X = (x,y)$ and $\bar{X} = (\bar{x}, \bar{y})$.

$$||\mathbb{F}(X) - \mathbb{F}(\bar{X})|| = |\mathbb{F}_1(X) - \mathbb{F}_1(\bar{X})| + |\mathbb{F}_2(X) - \mathbb{F}_2(\bar{X})|$$
$$= |\rho(x-\bar{x}) - \rho(x^2 - \bar{x}^2) - \varepsilon_1(x-\bar{x}) - (xy - \bar{x}\bar{y})|$$
$$+ \left|\frac{\psi xy(1+\phi\bar{x}) - \psi\bar{x}\bar{y}(1+\phi x)}{(1+\phi x)(1+\phi\bar{x})} - (1+\varepsilon_2)(y-\bar{y})\right|$$
$$\leq (\rho(1+M) + \varepsilon_1)|x-\bar{x}| + M|x-\bar{x}|$$
$$+ |\psi(xy(1+\phi\bar{x}) - \bar{x}\bar{y}(1+\phi x)) - (1+\varepsilon_2)(y-\bar{y})|$$
$$\leq (\rho + \varepsilon_1 + (\rho+1)M)|x-\bar{x}|$$
$$+ \psi|xy(1+\phi\bar{x}) - \bar{x}\bar{y}(1+\phi x)| + (1+\varepsilon_2)|y-\bar{y}|$$
$$\leq (\rho + \varepsilon_1 + (\rho+1)M)|x-\bar{x}|$$
$$+ \psi|xy - \bar{x}\bar{y}| + \psi\phi M|xy - \bar{x}\bar{y}| + (1+\varepsilon_2)|y-\bar{y}|$$
$$\leq (\rho + \varepsilon_1 + (\rho+1)M)|x-\bar{x}|$$
$$+ \psi M|y-\bar{y}| + \psi\phi M^2|y-\bar{y}| + (1+\varepsilon_2)|y-\bar{y}|$$
$$\leq L_1|x-\bar{x}| + L_2|y-\bar{y}| = L||X-\bar{X}|| \quad (7)$$

where $L_1 = \rho + \varepsilon_1 + (\rho+1)M$, $L_2 = 1 + \varepsilon_2 + \psi(1+\phi M)M$ and $L = \max\{L_1, L_2\}$. That is, with initial condition $X(t_0) = (x(t_0), y(t_0))$, an unique solution $X(t) \in S$ exists for system (4). □

**Theorem 3.2.** *Assume that $\Omega = \{(x,y,z) \in \mathbb{R}^3_+ : \max\{|x|,|y|,|z|\} \leq M\}$ and $S = \Omega \times [t_0, T]$ where $T < +\infty$. Then for any initial conditions $(x(t_0), y(t_0), z(t_0)) \in \Omega$, all the solutions $(x(t), y(t), z(t)) \in S$ of model (6) are non-negative and unique for all $t \geq 0$.*

**Proof.** (Proof by contradiction): let $k(t) = \min\{x(t), y(t), z(t)\}$, then $k(t) > 0$. Also assume

$$\exists \bar{t} > 0 \ s.t \ k(\bar{t}) = 0, \ \text{and} \ k(t) > 0 \ \forall t \in [0, \bar{t}).$$

If $k(\bar{t}) = x(\bar{t})$, then from the first equation of system (6), we have

$$x(\bar{t}) = x(0)E_\alpha\left[\left(\rho(1-x) - (y + \eta z + yz + \varepsilon_1)\right)\right] > 0$$

This leads to a contradiction. In the similar way, for $k(\bar{t}) = y(\bar{t})$ and $k(\bar{t}) = z(\bar{t})$, we can obtain the contradiction. Therefore, $\forall t \geq 0, \ k(t) > 0$ and as a result, $(x(t), y(t), z(t))$ will be positive for all $t \geq 0$.

As for the existence and uniqueness, after following the same procedure as Theorem 3.1 and considering a mapping function (Eq. (9)) as $\mathbb{F}(X) = (\mathbb{F}_1(X), \mathbb{F}_2(X), \mathbb{F}_3(X))$, we have found that under the following $L_i$ ($i = 1, 2, 3$), equation (10) holds for all $X = (x, y, z)$ and $\bar{X} = (\bar{x}, \bar{y}, \bar{z})$.

$$L_1 = \rho + \varepsilon_1 + (\rho + \eta)M + M^2,$$
$$L_2 = 1 + \varepsilon_2 + \psi\left[1 + (1+\phi)M + \phi M^2\right]M, \quad (8)$$
$$L_3 = \varepsilon_3 + \beta\left[1 + \eta + \eta\phi_1 M + \phi_1 M^2\right]M.$$

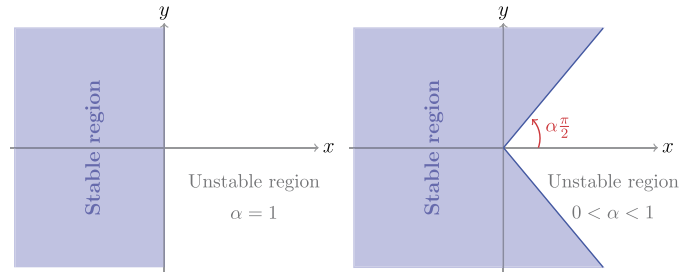

**Fig. 1.** Stability region of the fractional-order system.

$$\mathbb{F}_1(X) = \rho x(1-x) - x(y + \eta z + yz) - \varepsilon_1 x,$$
$$\mathbb{F}_2(X) = \frac{\psi xy(1+z)}{1+\phi x} - \varepsilon_2 y, \quad (9)$$
$$\mathbb{F}_3(X) = \frac{\beta xz(\eta + y)}{1+\phi_1 x} - \varepsilon_3 z.$$

$$||\mathbb{F}(X) - \mathbb{F}(\bar{X})|| \leq L||X - \bar{X}|| \ \text{s.t} \ L = \max\{L_1, L_2, L_3\} \quad (10)$$

Hence, with initial condition $X(t_0) = (x(t_0), y(t_0), z(t_0))$, an unique solution $X(t) \in S$ exists for system (6). □

### 3.2. Stability analysis

In this section, we have examined the stability of equilibrium solutions of the proposed models. All the calculated stability conditions are summarised into Table 1 and Table 2.

To obtain the stability of dynamic models with integer order, one very common way is to take advantage of the Routh-Hurwitz (RH) conditions.

i) If $\det(\mathbf{J}(E^\star)) > 0$, then $E^\star$ will be a saddle point (regardless of the sign of $\text{Tr}(\mathbf{J})$),
ii) If $\det(\mathbf{J}(E^\star)) < 0$, then $E^\star$ will be asymptotically stable when $\text{Tr}(\mathbf{J}(E^\star)) < 0$,
iii) If $\det(\mathbf{J}(E^\star)) < 0$, then $E^\star$ is unstable when $\text{Tr}(\mathbf{J}(E^\star)) > 0$.

where $\mathbf{J}$ and $E^\star$ are the Jacobian matrix and equilibrium point of the dynamic system respectively.

In addition to RH conditions, in fractional calculus, one needs to examine another condition (Lemma 3.1) [41].

**Lemma 3.1** ([37, page 158]). *Consider function $D^\alpha F(x) = g(x)$ with $0 < \alpha \leq 1$. An equilibrium point ($E^\star$) is asymptotically stable if all eigenvalues ($\lambda_i, i = 1, 2, ..., n$) of the associated Jacobian matrix in $E^\star$ satisfy the following condition (Fig. 1)*

$$|\arg(\lambda_i)| \geq \frac{\alpha\pi}{2}, \quad i = 1, 2, ..., n. \quad (11)$$

What stands out from the Fig. 1 is that the stability domain of the fractional order system is larger than the corresponding domain for integer systems when $\alpha \in (0, 1)$. Therefore it is expected that a non-integer derivative will increase the stability of the system.

**Definition 3.2.** Consider the following polynomial

$$f(x) = x^n + a_1 x^{n-1} + a_2 x^{n-2} + \cdots + a_n \quad (12)$$

the discriminant of the polynomial $f(x)$ is defined by $D(f) = (-1)^{n(n-1)/2} R(f, f')$ where $f'$ is the derivative of $f$ and where $g(x) = x^n + b_1 x^{l-1} + b_2 x^{l-2} + \cdots + b_l$ and $R(f, g)$ is an $(n+l) \otimes (n+l)$ determinant.

For $n = 3$, $D(f) > 0$ implies that all the roots are real and $D(f) < 0$ implies that there is only one real root and one complex and its complex conjugate. Also

$$D(f) = 18a_1 a_2 a_3 + (a_1 a_2)^2 - 4a_3(a_1)^3 - 4(a_2)^3 - 27(a_3)^2.$$

**Lemma 3.2** ([41]). *Consider characteristic equation (12), then the conditions that satisfies (11) are as follows.*
i) For $n = 1$, the condition is $a_1 > 0$.
ii) For $n = 2$, the conditions are either RH conditions or

$$a_1 < 0, \quad 4a_2 > (a_1)^2, \quad |\tan^{-1}(\sqrt{4a_2 - (a_1)^2})/a_1| > \alpha \pi/2.$$

iii) For $n = 3$, if $D(f) < 0$, $a_1 < 0, a_2 < 0, a_3 > 0$, then (11) is satisfied for all $\alpha > (2/3)$.
iv) If $D(f) > 0$, $a_1 < 0, a_2 < 0, a_3 > 0$, then (11) is satisfied for all $0 < \alpha < 1$.
v) For general $n$, $a_n > 0$, is a necessary condition for (11).

*3.2.1. Single predator model*

The system (4) has three stationary points: $E_1 = (0,0)$ extinction of both species; $E_2 = (1 - \varepsilon_1/\rho, 0)$ predator only extinction; and $E_3 = (\omega, \rho(1 - \omega) - \varepsilon_1)$ predator/prey co-existence, where

$$\omega = \frac{1 + \varepsilon_2}{\psi - \phi(1 + \varepsilon_2)}.$$

According to the stationary points, non-zero harvesting clearly shifts the prey-only equilibrium, $E_2$, to lower densities of prey. For the co-existence point, $E_3$, the prey population density is unaffected by harvesting of the prey itself but is shifted to higher densities as the rate of predator harvesting is increased. The harvesting of prey is effectively competition for the predators reducing their equilibrium population in concert with any direct harvesting. Also

$$\mathbf{J}(x^\star, y^\star) = \begin{pmatrix} \rho(1 - 2x^\star) - y^\star - \varepsilon_1 & -x^\star \\ \frac{\psi y^\star}{(1 + \phi x^\star)^2} & \frac{(\psi - \phi)x^\star - 1}{1 + \phi x^\star} - \varepsilon_2 \end{pmatrix}. \quad (13)$$

where $\mathbf{J}(x^\star, y^\star)$ is the Jacobian matrix of the system (4) in the stationary points $(x^\star, y^\star)$.

From (13), it can easily be deduced that $\lambda_1 = \rho - \varepsilon_1$ and $\lambda_2 = -(1 + \varepsilon_2)$ in extinction point ($E_1$). According to Lemma 3.1 therefore, $E_1$ is asymptotically stable for $0 < \alpha \le 1$ if and only if $\rho < \varepsilon_1$. That is, if the rate of harvesting of the prey outstrips its growth rate, the system is driven to extinction.

Also, it can be concluded that the linearised system has two eigenvalues of

$$\lambda_1 = -\rho\left(1 + \frac{\varepsilon_1}{\rho}\right) - \varepsilon_1 \quad \text{and} \quad \lambda_2 = \frac{\psi(\rho - \varepsilon_1)}{\rho + \phi(\rho - \varepsilon_1)} - (1 + \varepsilon_2)$$

at $E_2$, and thus according to Lemma 3.1, $E_2$ is asymptotically stable for all $0 < \alpha \le 1$ only when

$$\psi < (1 + \varepsilon_2)\left[\frac{\rho}{(\rho - \varepsilon_1)} + \phi\right],$$

which corresponds to the coexistence point being biologically irrelevant. For fixed point $E_3$, we have

$$\mathbf{J}(E_3) = \begin{pmatrix} -(\rho\omega + \varepsilon_1) & -\omega \\ \frac{\psi(\rho(1 - \omega) - \varepsilon_1)}{(1 + \phi\omega)^2} & 0 \end{pmatrix}. \quad (14)$$

Thus, the trace and determinant of the Jacobian evaluated at Eq. (14) are as follows:

$$\text{Tr}(\mathbf{J}(E_3)) = -(\rho\omega + \varepsilon_1), \quad \det(\mathbf{J}(E_3)) = \frac{\psi\omega(\rho(1 - \omega) - \varepsilon_1)}{(1 + \phi\omega)^2}.$$

According to RH, the sufficient condition for asymptotically stability of the system (4) is when $\omega < 1 - \varepsilon_1/\rho$.

To satisfy the Lemma 3.1 as the necessary condition, we must have (Fig. 2)

$$((\rho + \varepsilon_1 \phi)(1 + \varepsilon_2) - \varepsilon_1 \psi)^2 < \frac{4}{\psi} \rho(1 + \varepsilon_2)(\psi - \phi(1 + \varepsilon_2))^2$$
$$\times \left[(1 - \varepsilon_1)(\psi - \phi(1 + \varepsilon_2)) - (1 + \varepsilon_2)\right]. \quad (15)$$

*3.3. Two-predator model*

There are five stationary points for model (6):

- $E_1 = (0,0,0)$, $E_2 = (1 - \frac{\varepsilon_1}{\rho}, 0, 0)$;
- $E_3 = (x^\star, \rho(1 - x^\star) - \varepsilon_1, 0)$ s.t $x^\star = \frac{\varepsilon_2}{\psi - \varepsilon_2 \phi}$;
- $E_4 = (x^\star, 0, \frac{1}{\eta}(\rho(1 - x^\star) - \varepsilon_1))$ s.t $x^\star = \frac{\varepsilon_3}{\eta\beta - \phi_1 \varepsilon_3}$; and
- $E_5 = (x^\star, \frac{\varepsilon_3(1 + \phi_1 x^\star)}{\beta x^\star} - \eta, \frac{\varepsilon_2(1 + \phi x^\star)}{\psi x^\star} - 1)$ s.t $x^\star = 1 - \frac{1}{\rho}\left[\frac{\varepsilon_2 \varepsilon_3 (1 + \phi\omega)(1 + \phi_1 \omega)}{\beta\psi(\omega)^2} + \varepsilon_1 - \eta\right]$.

where $E_1$ is total population extinction, $E_2$ is prey only, $E_3$ and $E_4$ are equilibria of partial co-existence (the prey with one of the predators), and $E_5$ is the co-existence of all three species with $\omega \in \mathbb{R}^+$.

Linearising the system (6) about the stationary points $(x^\star, y^\star, z^\star)$, we can determine each point's linear stability by considering the eigenvalues of the resulting Jacobian matrix.

$$\mathbf{J}(x^\star, y^\star, z^\star)$$
$$= \begin{pmatrix} \rho(1 - 2x^\star) - (y^\star(1 + z^\star) + \eta z^\star + \varepsilon_1) & -x^\star(1 + z^\star) & -x^\star(\eta + y^\star) \\ \frac{\psi y^*(1 + z^\star)}{(1 + \phi x^\star)^2} & \frac{\psi x^\star(1 + z^\star)}{1 + \phi x^\star} - \varepsilon_2 & \frac{\psi x^\star y^\star}{1 + \phi x^\star} \\ \frac{\beta z^\star(\eta + y^\star)}{(1 + \phi_1 x^\star)^2} & \frac{\beta x^\star z^\star}{1 + \phi_1 x^\star} & \frac{\beta x^\star(\eta + y^\star)}{1 + \phi_1 x^\star} - \varepsilon_3 \end{pmatrix}.$$
(16)

From (16) and Lemma 3.1, it can easily be deduced that $E_1$ is asymptotically stable for $0 < \alpha \le 1$ if only if $\rho < \varepsilon_1$. At $E_2$ also, we have three following eigenvalues

$$\lambda_1 = -\rho + \varepsilon_1, \quad \lambda_2 = \frac{\psi x^\star}{1 + \phi x^\star} - \varepsilon_2, \text{ and} \quad \lambda_3 = \frac{\beta\eta x^\star}{1 + \phi_1 x^\star} - \varepsilon_3$$

thus according to Lemma 3.1, $E_2$ is asymptotically stable for all $0 < \alpha \le 1$ only when

$$\varepsilon_1 < \rho, \quad \varepsilon_2 > \frac{\psi(\rho - \varepsilon_1)}{\rho + \phi(\rho - \varepsilon_1)}, \text{ and} \quad \varepsilon_3 > \frac{\beta\eta(\rho - \varepsilon_1)}{\rho + \phi_1(\rho - \varepsilon_1)}$$

As for $E_3$, we have the following characteristic equation.

$$\det(\mathbf{J}(E_3) - \lambda I) = \left[\frac{\beta x^\star(\eta + y^\star)}{1 + \phi_1 x^\star} - (\varepsilon_3 + \lambda)\right] \begin{vmatrix} \rho x^\star - \lambda & -x^\star \\ \frac{\psi y^\star}{(1 + \phi x^\star)^2} & -\lambda \end{vmatrix}.$$

Therefore, according to the Lemma 3.1 and RH, the necessary and sufficient conditions for asymptotically stability of the system (6) at $E_3$ for all $\alpha \in (0, \alpha_1)$ is

$$\varepsilon_1 < \rho(1 - \omega), \quad \varepsilon_2 < \frac{\psi}{\phi}, \quad \varepsilon_3 > \frac{\beta\omega(\eta + \gamma)}{1 + \phi_1 \omega}, \text{ and} \quad \rho < \frac{2}{(1 + \phi\omega)}\sqrt{\frac{\gamma}{\omega}\psi}$$

where $\omega = \frac{\varepsilon_2}{\psi - \varepsilon_2 \phi}$, $\gamma = \rho(1 - \omega) - \varepsilon_1$ and $\alpha_1 = \min\{1, \frac{2}{\pi}|\arg(-\frac{\rho\omega}{2} \pm i\frac{\sqrt{4\psi\omega\gamma/(1+\phi\omega)^2 - (\rho\omega)^2}}{2})|\}$.

Regarding $E_4$, according to (16), we will also have the following characteristic equation.

$$\det(\mathbf{J}(E_4) - \lambda I) = \left[\frac{\psi x^\star(1 + z^\star)}{1 + \phi x^\star} - (\varepsilon_2 + \lambda)\right] \begin{vmatrix} \rho x^\star - \lambda & -\eta x^\star \\ \frac{\beta\eta z^\star}{(1 + \phi_1 x^\star)^2} & -\lambda \end{vmatrix}.$$

Using the RH condition (ii) and the Lemma 3.1, $E4$ is asymptotically stable for all $\alpha \in (0, \alpha_1)$ only if

$$\varepsilon_1 < \rho(1 - \omega), \quad \varepsilon_2 > \frac{\psi\omega(\gamma + 1)}{\eta(1 + \phi\omega)}, \quad \varepsilon_3 < \frac{\eta}{\phi_1}\beta, \text{ and} \quad \rho < \frac{2}{(1 + \phi_1 \omega)}\sqrt{\frac{\gamma}{\omega}\eta\beta}$$

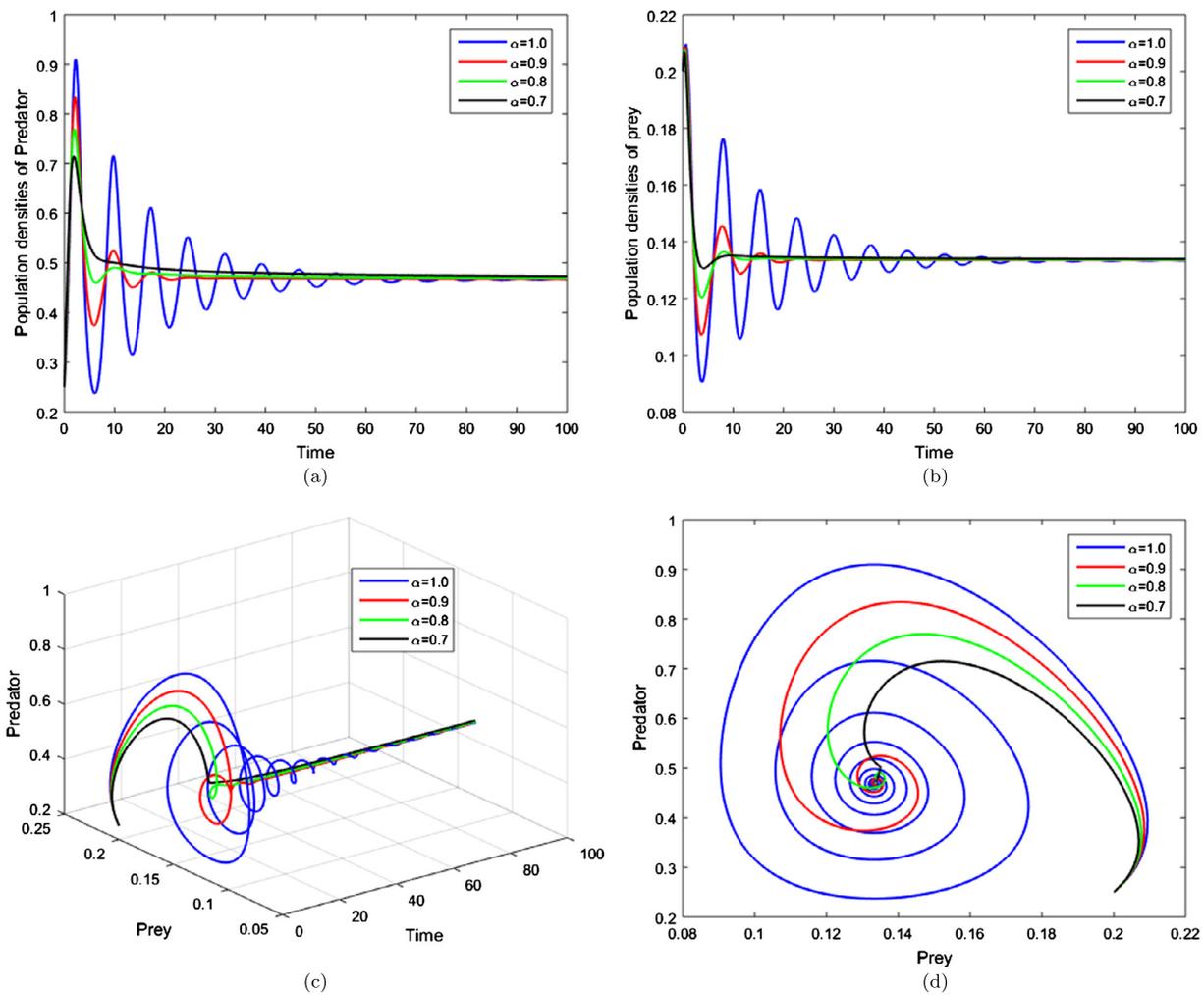

**Fig. 2.** The damping of population oscillations produced by system (4) with initial densities (0.2,0.25) and $\psi = 19$, $\rho = 1$, $\psi = 15$, $\phi = 2$, $\varepsilon_1 = 0.4$ and $\varepsilon_2 = 1$. (a and b) Numerical values of the lanternfish and tuna versus time respectively. (c): Combined predator-prey population densities with respect to time and (d): Phase portrait of the system.

where $\omega = \frac{\varepsilon_3}{\eta\beta - \phi_1\varepsilon_3}$, $\gamma = \rho(1-\omega) - \varepsilon_1$, and $\alpha_1 = \min\left\{1, \frac{2}{\pi}\left|\arg(-\frac{\rho\omega}{2} \pm i\frac{\sqrt{4\eta\beta\omega\gamma/(1+\phi_1\omega)^2 - (\rho\omega)^2}}{2})\right|\right\}$.

Also, the following characteristic equation is given for $E_5$.

$$\det\left(\mathbf{J}(E_5) - \lambda I\right)$$
$$= \lambda^3 - \left[y^\star(1+z^\star) + \eta z^\star + \rho\right]\lambda^2 + \left[\frac{\varepsilon_3^2\psi z^\star + \varepsilon_2^2\beta y^\star}{\beta\psi x^\star}\right.$$
$$\left. - \frac{\psi\beta(x^\star)^2 y^\star z^\star}{(1+\phi_1 x^\star)(1+\phi x^\star)}\right]\lambda$$
$$+ \left[\varepsilon_2\varepsilon_3 y^\star z^\star + \frac{\rho\psi\beta(x^\star)^2 y^\star z^\star + \varepsilon_2\varepsilon_3 y^\star z^\star(2+(\phi+\phi_1)x^\star)}{(1+\phi x^\star)(1+\phi_1 x^\star)}\right],$$

thus according to Lemma 3.2, $E_5$ under the following conditions are asymptotically stable for all $\alpha \in (0,1)$.

$\varepsilon_1 < \rho + \eta - \gamma\gamma_1$, $\varepsilon_2 > \frac{\psi\omega_1}{1+\phi\omega_1}$, $\varepsilon_3 > \frac{\eta\beta\omega_1}{1+\phi_1\omega_1}$, and

$\psi\varepsilon_3^2 + \eta\beta\varepsilon_2^2 > \psi\varepsilon_3^2\gamma_1 + \beta\varepsilon_2^2\gamma - \frac{\psi\omega}{\gamma\gamma_1}\beta\varepsilon_1\varepsilon_2$,

where $\gamma = \frac{\varepsilon_3(1+\phi_1\omega)}{\beta\omega}$, $\gamma_1 = \frac{\varepsilon_2(1+\phi\omega)}{\psi\omega}$, $\omega_1 = 1 - \frac{1}{\rho}\left[\gamma\gamma_1 + \varepsilon_1 - \eta\right]$, and $\omega \in \mathbb{R}^+$ (Fig. 3d).

## 4. Physical interpretation of the models

In this section, we have numerically [42] simulated both systems (4) and (6) under the stability conditions found in the previous section. Also, we presented some possible physical explanations for the obtained results.

From the simulations presented in Figs. 2 and 3, the effect of reducing the order of the time derivative can be seen. As the fractional order $\alpha$ is decreased, the system (with Caputo derivative) stabilises faster. That is the higher "memory" the system has of past states, the greater the damping of oscillations in system dynamics. The simulations demonstrate that, even with quite moderate reductions in $\alpha$, the amplitude of population density oscillations is strongly retarded.

From the numerical simulations we see that the harvesting appears to enhance the stability of the system. When no harvesting is conducted, oscillations in the population numbers are immense (Fig. 4a). Once at least one of the harvesting parameters is non-zero these oscillations are damped out (Fig. 4b), aligning with the earlier results presented in the literature [15, 20].

We emphasize the novel relationship between memory concept and harvesting. The numerical integration of system (4) (blue lines) indicates that when interaction among species is immense, such that the population numbers fluctuate wildly, harvesting could be a solution to stabilise the system. For some parameter values, harvesting works less well as a stabilising mechanism and could destabilise the coexistence fixed point of the system (Fig. 4c). In other words, the modelling assumption of the "largeness" of the population, allowing the use of differential equations, may very well breakdown. In these situations, harvesting reduces the range of values through which the population numbers fluctuate (thus reducing the time for the system to settle at equilibrium) but also shifts the range to lower values. The stability of the fixed point as per the stability condition is likely incorrect because the population returns to healthy levels after having dropped to number better interpreted as having gone extinct (Fig. 4c). To smooth out this inconsistency, our findings suggest that another parameter, defined as species' memory, should be incorporated in the model (Fig. 4d). A possible explanation for this could that when species are heavily harvested, they rely much more on their memory to prevent the species extinction. For example, if they are heavily harvested during the spring season, they produce more offspring intrinsically during that period to make up their lack of population. In other words, they might have a variable growth rate to compensate for the different physical circumstances from one breeding season to another; they might find a safer place to reproduce so that they reduce the death-rate of infants, or even migrate to another place where they can decrease mortality rates of offspring through less exposure to predators. Therefore, we believe that such scenarios could be explained by the memory concept.

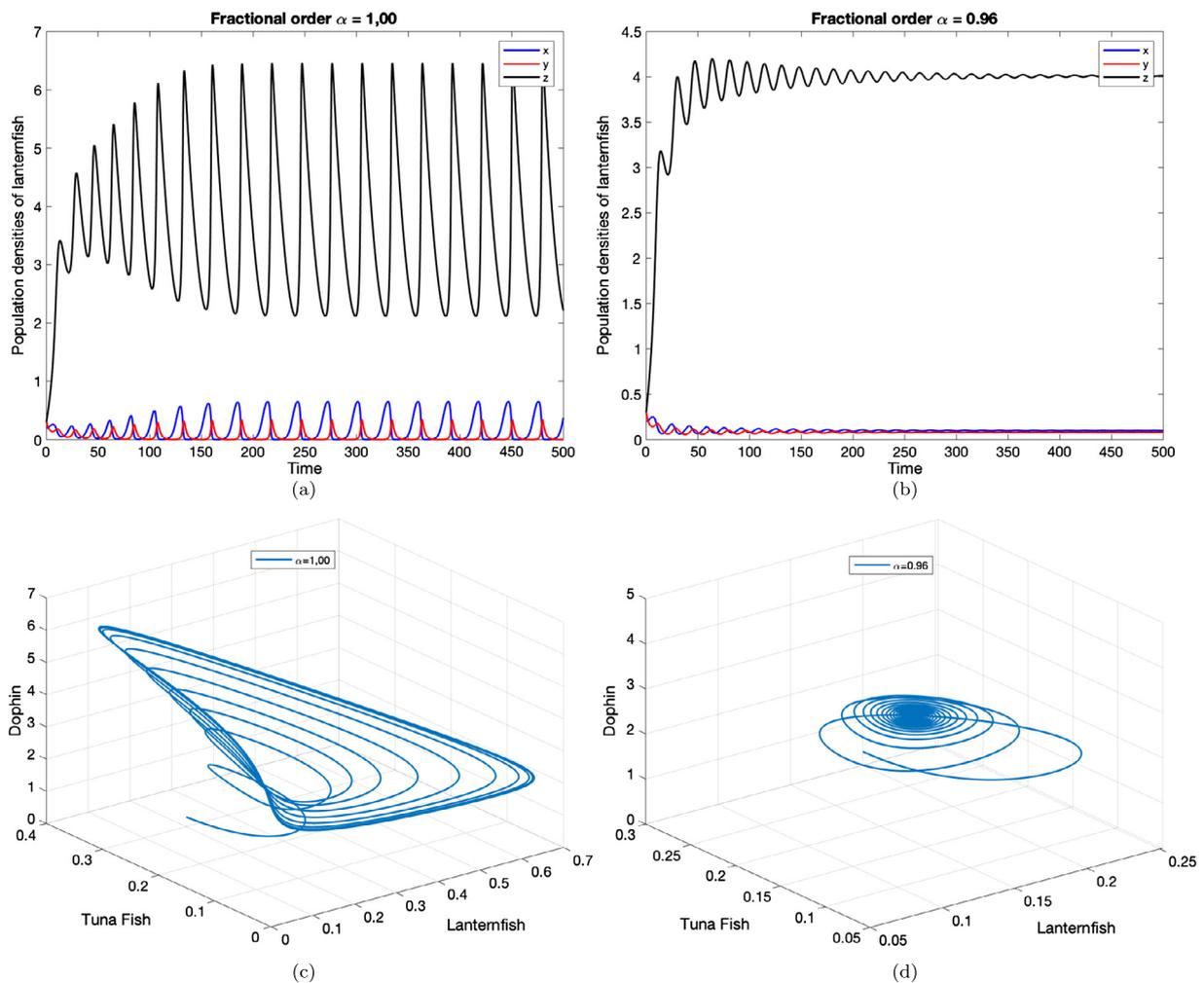

**Fig. 3.** The damping of population oscillations produced by system (6) with integer (a) and fractional order (b). (c and d): Phase portrait of the integer and fractional system respectively. The initial densities $(0.2, 0.3, 0.3)$, $(\rho, \psi, \beta) = (0.61, 1, 7)$, $(\phi, \phi_1) = (1.4, 0.02)$, $(\varepsilon_1, \varepsilon_2, \varepsilon_3) = (0.12, 0.43, 0.06)$ and $\eta = 0.01$.

## 5. Discussion

The mutualistic interactions amongst species is of great importance in the field of conservation ecology, so gaining an understanding of such interactions can make a noteworthy contribution to species maintenance. With this in mind, we introduced a modified Lotka-Volterra model to study interaction amongst three species with mutualistic predation. The motivation is based on the observed feeding behaviour of spotted dolphins and yellowfin tuna upon schools of lanternfish. Along with obtaining the stability conditions for the model, we also investigated the impact of "memory effects" on the species interaction via fractional calculus. Our analysis reveals that the fractional system dampens out induced oscillating inherent in predator-prey models and reaches the local stable point sooner than the integer model does. In other words, stability is more robust when the species exhibit "memory".

Moreover, assuming that either both species have market value or that one species is caught as by-catch, we have investigated the effect of constant rate harvesting within the proposed models. We have also discussed the local stability behaviour of all the equilibrium states of the system. The output of all models shows that stability and extinction of the ecosystem are affected by economic interest/harvesting. In fact, applying constant harvest at a rate below the threshold of the stability condition causes the system to stabilise faster. In other words, exploitation of a species can be regarded as beneficial to the ecosystem as a whole as the system will reach steady-state sooner. Further, the greater the collaboration among species (large mutualism coefficients), the more harvesting that can occur before putting the populations at risk of extinction (over-fished).

An assumption of a differential equation based population models is that the quantities under consideration are sufficiently large such that the addition or loss of an individual to that population may be considered an infinitesimal change. That is, the population may take on a continuum of values rather than being strictly discrete. Further, it is well known that many such models can exhibit, for certain parameter values and initial conditions, large fluctuations — particularly those with Hopf bifurcations — whereby the assumption of the "largeness" of the population may very well breakdown Fig. 5. In practical terms the population has either already gone extinct, even though the model predicts a return to healthy levels, or its reduced to a level vulnerable to external perturbations such that it future viability should be discussed in probabilistic terms. These issues imply the model should incorporate stochasticity at low population levels. Alternatively, assuming a memory like behaviour can be attributed to the population, a non-integer order time derivative arguably extends the validity of the model by preventing the wild swings in population numbers that represent 'essentially extinct' to 'thriving' (Fig. 5a).

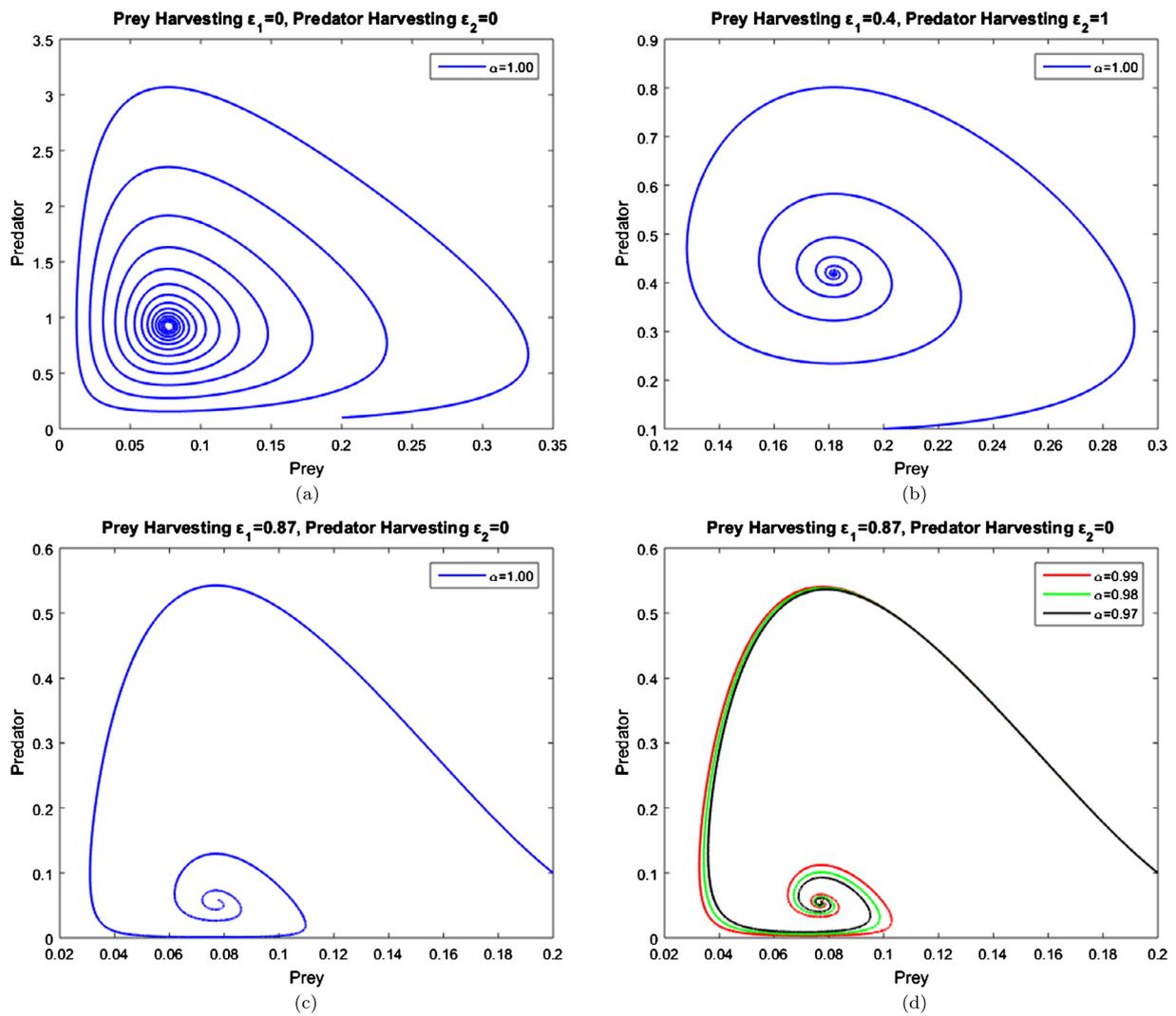

**Fig. 4.** The impact of harvesting upon the system (4) with fractional order ($0 < \alpha \leq 1$), $\rho = 1, \psi = 15, \phi = 2$ and $(x(0), y(0)) = (0.2, 0.1)$.

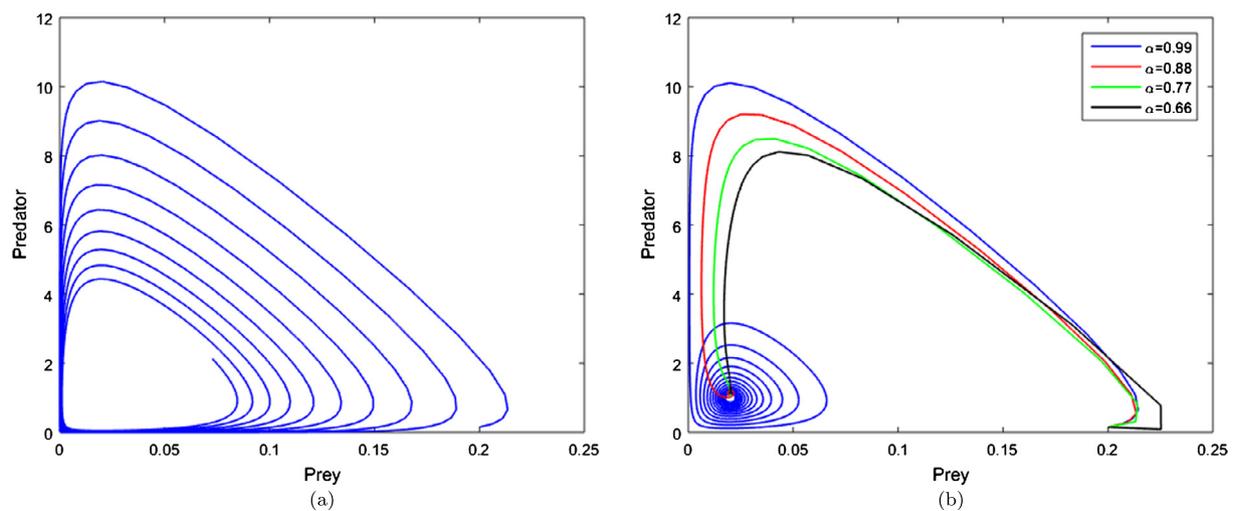

**Fig. 5.** Phase plane of system (4) with integer order derivative (5a) and fractional order derivative (5b). All initial densities are $(x, y) = (0.2, 0.15)$ with parameter values $(\rho, \psi, \phi, \varepsilon_1, \varepsilon_2) = (1, 50, 0.02, 0, 0)$. For these parameter values criterion (15) holds.

It has been argued that ecological systems, while being stable for moderate numbers of species interactions or moderate strength of the inter-species connections, will become unstable to small perturbations once some threshold value for interaction number or interaction strength is breached [43, 44]. In 2012, similar findings were presented, claiming that any feasible co-existence of the system species will be unstable when the pair-wise competitive interactions are sufficiently strong [45]. Further, Goh presented findings that the continuum of globally stable Lotka-Volterra models of mutualism among three or more species is smaller than the continuum of globally stable Lotka-Volterra models of competition among the same number of species [46]. Thus, if diversity has an adverse effect on stability in competitive system it has even more so in mutualistic ones.

More recently, however, Butler and O'Dwyer have to some extent overturned this understanding through proposing a consumer-producer-resource model. They demonstrated, for a model of $N$ bacteria species consuming $N$ abiotic resources that the stability is guaranteed for all feasible equilibria. For an extension of the model where the bacteria also produce some or all of the resources to mutual benefit of all consumers, stability of all feasible solutions can be guaranteed provided that mutualistic interactions are symmetric [47].

Our system, on the other hand, is one with moderate to strong interactions between mutualistic predators and a prey (biotic resource) species. Our formulation allows weak to strong memory effects via the fractional derivative order. The three-species model suggests the stability at the co-existence is not guaranteed by symmetric mutualism. Rather, the model needs to have a memory to satisfy the stability condition. As future work we will investigate how stability of feasible solutions in the model changes as the number of predators and prey increase.

## Declarations

### Author contribution statement

M.M. Amirian: Conceived and designed the experiments; Performed the experiments; Analyzed and interpreted the data; Wrote the paper.

I.N. Towers, Z. Jovanoski, A.J. Irwin: Analyzed and interpreted the data; Wrote the paper.

### Funding statement

This research did not receive any specific grant from funding agencies in the public, commercial, or not-for-profit sectors.

### Competing interest statement

The authors declare no conflict of interest.

**Table 1.** Asymptotic stability analysis of system (4). The notation ✓ indicates that the relevant condition is satisfied.

| | System (4) | |
|---|---|---|
| $E_i$: Equilibrium point | Stability conditions | Fig. 2 |
| $E_1$: Extinction point | $\rho < \varepsilon_1$ | … |
| | $\alpha \leq 1$ | |
| $E_2$: Predator-free point | $\psi < (1+\varepsilon_2)\left[\frac{\rho}{(\rho-\varepsilon_1)}+\phi\right]$ | … |
| | $\alpha \leq 1$ | |
| $E_3$: Coexistence points | $\omega < 1 - \frac{\varepsilon_1}{\rho}$ | $0.13 < 1 - 0.4 \rightarrow$ ✓ |
| | $b^2 - 4c < 0$ | $(0.53)^2 - 3.7 < 0 \rightarrow$ ✓ |
| | $\alpha < \frac{2}{\pi}\lvert\arg(-\frac{b}{2} \pm i\frac{\sqrt{4c-b^2}}{2})\rvert$ | $\alpha \in (0, 1.178) \rightarrow$ ✓ |

### Additional information

No additional information is available for this paper.

## Acknowledgements

MA, the first author, expresses his appreciation for many useful and enlightening discussions with *Dr. Yousef Jamali* of Tarbiat Modares University and *Moein Khalighi* of the University of Turku.

## Appendix A

The stability conditions given in section 3, are summarised in two general tables. In Table 1, we have defined the following notations:

$$b = -\text{Tr}(\mathbf{J}(E_3)), c = \det(\mathbf{J}(E_3)), \text{ and } \omega = (1+\varepsilon_2)/(\psi - \phi(1+\varepsilon_2)).$$

In Table 2 also, we have the following notations for $E_i (i = 3, 4, 5)$:

- In $E_3$, known as partial co-existence point

$$\omega = \frac{\varepsilon_2}{\psi - \varepsilon_2 \phi}, \gamma = \rho(1-\omega) - \varepsilon_1, \text{ and}$$

$$\alpha_1 = \min\{1, \frac{2}{\pi}\lvert\arg(-\frac{\rho\omega}{2} \pm i\frac{\sqrt{4\psi\omega\gamma/(1+\phi\omega)^2 - (\rho\omega)^2}}{2})\rvert\}$$

- In $E_4$, known as partial co-existence point

$$\omega = \frac{\varepsilon_3}{\eta\beta - \phi_1\varepsilon_3}, \gamma = \rho(1-\omega) - \varepsilon_1, \text{ and}$$

$$\alpha_1 = \min\{1, \frac{2}{\pi}\lvert\arg(-\frac{\rho\omega}{2} \pm i\frac{\sqrt{4\eta\beta\omega\gamma/(1+\phi_1\omega)^2 - (\rho\omega)^2}}{2})\rvert\}$$

- In $E_5$, known as co-existence point

$$\omega_1 = 1 - \frac{1}{\rho}\left[\gamma\gamma_1 + \varepsilon_1 - \eta\right], \omega \in \mathbb{R}^+, \gamma = \frac{\varepsilon_3(1+\phi_1\omega)}{\beta\omega}, \text{ and } \gamma_1 = \frac{\varepsilon_2(1+\phi\omega)}{\psi\omega}.$$

**Table 2.** Asymptotic stability analysis of system (6). The notation ✓ indicates that the relevant condition is satisfied.

| $E_i$: Equilibrium point | Stability conditions | Fig. 3 |
|---|---|---|
| **System (6)** | | |
| $E_1$: **Extinction point** $(0,0,0)$ | $\rho < \varepsilon_1$ <br> $\alpha \leq 1$ | … |
| $E_2$: **Predator-free point** $(x^\star, 0, 0)$ | $\varepsilon_1 < \rho$ <br> $\varepsilon_2 > \psi(\rho - \varepsilon_1)/[\rho + \phi(\rho - \varepsilon_1)]$ <br> $\varepsilon_3 > \beta\eta(\rho - \varepsilon_1)/[\rho + \phi_1(\rho - \varepsilon_1)]$ <br> $\alpha \leq 1$ | … |
| $E_3$: **Partial coexistence point** $(x^\star, y^\star, 0)$ | $\varepsilon_1 < \rho(1-\omega)$ <br> $\varepsilon_2 < \psi/\phi$ <br> $\varepsilon_3 > \beta\omega(\eta+\gamma)/[1+\phi_1\omega]$ <br> $\rho < \frac{2}{(1+\phi\omega)}\sqrt{\frac{\gamma}{\omega}\psi}$ <br> $\alpha \in (0, \alpha_1)$ | … |
| $E_4$: **Partial coexistence point** $(x^\star, 0, z^\star)$ | $\varepsilon_1 < \rho(1-\omega)$ <br> $\varepsilon_2 > \psi\omega(\gamma+1)/[\eta(1+\phi\omega)]$ <br> $\varepsilon_3 < \eta\beta/\phi_1$ <br> $\rho < \frac{2}{(1+\phi_1\omega)}\sqrt{\frac{\gamma}{\omega}\eta\beta}$ <br> $\alpha \in (0, \alpha_1)$ | … |
| $E_5$: **Coexistence point** $(x^\star, y^\star, z^\star)$ | $\varepsilon_1 < \rho + \eta - \gamma\gamma_1$ <br> $\varepsilon_2 > \psi\omega_1/[1+\phi\omega_1]$ <br> $\varepsilon_3 > \eta\beta\omega_1/[1+\phi_1\omega_1]$ <br> $\psi\varepsilon_3^2 + \eta\beta\varepsilon_2^2 > \psi\varepsilon_3^2\gamma_1 + \beta\varepsilon_2^2\gamma - \frac{\psi\omega}{\gamma\gamma_1}\beta\varepsilon_1\varepsilon_2$ <br> $\alpha \in (0,1)$ | $0.12 < 0.62 \to$ ✓ <br> $0.43 > 0.3812 \to$ ✓ <br> $0.06 > 0.0563 \to$ ✓ <br> $0.016 > -879.9 \to$ ✓ <br> $0.96 \in (0,1) \to$ ✓ |